\documentclass[11pt,a4paper]{article}
\usepackage[cp1251]{inputenc}
\usepackage{amsfonts}
\usepackage{amssymb}
\usepackage{amsthm}
\usepackage{amsmath}
\usepackage{mathtext}
\usepackage{longtable}
\usepackage[russian]{babel}
\usepackage{color}

\theoremstyle{remark}

\newcommand{\R}{\mathbb R}
\renewcommand{\r}{\mathbb R}
\newcommand{\N}{\mathbb N}

\renewcommand{\le}{\leqslant}
\renewcommand{\ge}{\geqslant}

\renewcommand{\phi}{\varphi}

\newcommand{\eqd}{\stackrel{d}{=}}

\newcommand{\vp}{\emph{\textbf{p}}}
\newcommand{\wX}{\widetilde X}
\definecolor{lgreen}{rgb}{0.9,1,0.8}

\title{Bounds for the concentration functions of random sums under
relaxed moment conditions\thanks{Research supported by the Russian
Foundation for Basic Research, project 15-07-02984.}}
\author{V. Yu. Korolev\thanks{Hangzhou Dianzi University; Faculty of
Computational Mathematics and Cybernetics, Lomonosov Moscow State
University; Federal Research Center ``Informatics and Control'',
Russian Academy of Sciences; vkorolev@cs.msu.su}, A. V.
Dorofeeva\thanks{Faculty of Computational Mathematics and
Cybernetics, Lomonosov Moscow State University;
alex.dorofeyeva@gmail.com}}
\date{}

\begin{document}

\maketitle

{\bf Abstract.} Estimates are constructed for the deviation of the
concentration functions of sums of independent random variables with
finite variances from the folded normal distribution function
without any assumptions concerning the existence of the moments of
summands of higher orders. The obtained results are extended to
Poisson-binomial, binomial and Poisson random sums. Under the same
assumptions, the bounds are obtained for the approximation of the
concentration functions of mixed Poisson random sums by the
corresponding limit distributions. In particular, bounds are
obtained for the accuracy of approximation of the concentration
functions of geometric, negative binomial and Sichel random sums by
the exponential, the folded variance gamma and the folded Student
distribution. Numerical estimates of all the constants involved are
written out explicitly.

\smallskip

{\bf Key words:} distribution function, central limit theorem,
normal distribution, folded normal distribution, uniform metric,
Poisson-binomial distribution, Poisson-binomial random sum, binomial
random sum, Poisson random sum, mixed Poisson random sum, geometric
random sum, gamma distribution, negative binomial random sum,
inverse gamma distribution, Sichel distribution, Laplace
distribution, exponential distribution, folded variance gamma
distribution, folded Student distribution, absolute constant.

\normalsize



\renewcommand{\thesection}{}

\section{Introduction}

Assume that all the random variables considered in this paper are
defined on one and the same probability space
$(\Omega,\,\mathfrak{A},\,{\sf P})$.

The {\it concentration function} $Q_{\xi}(z)$ of a random variable
$\xi$ is defined as
$$
Q_{\xi}(z)=\sup_{x\in\r}{\sf P}(x\le\xi\le x+z),\ \ \ \ z\ge
0.\eqno(1)
$$
This notion was introduced by P. L{\'e}vy in 1937 \cite{Levy1937},
for details see \cite{Hengartner, Petrov1987}. Concentration
functions are convenient and informative characteristics of the
dispersion or scatter of random variables. It is conventional and
convenient to characterize the dispersion of a random variable by
its variance which is very simple to understand since it is a single
number. However, the cost of the simplicity of the variance is the
absence of the information concerning what deviations of the random
variable from its expected value are more probable than the others.
This information is contained in the concentration functions.

The estimates of the rate of decrease of the concentration functions
of sums of independent random variables as the number of summands
grows are well known, see, e. g., \cite{Petrov1987}. However, these
estimates are rather rough and do not take into consideration the
corresponding change of the shape of the concentration function.
Perhaps, for the first time the estimates of the concentration
functions that describe the asymptotic change of their shapes were
obtained in the paper \cite{Galieva2}. In these papers the estimates
were obtained under the assumption of existence of the third
moments.

In the present paper we relax the moment conditions and construct
the estimates of the concentration functions only under the
condition of existence of the variances of summands. The resulting
estimates make it possible to directly compare the informativeness
of the concentration function with that of the variance as the
measure of dispersion.

Along with purely theoretical motivation, there is a somewhat
practical interest in the problems considered below.
Poisson-binomial, binomial and mixed Poisson (first of all,
geometric) random sums are widely used as stopped-random-walk models
in many fields such as financial mathematics (Cox--Ross--Rubinstein
binomial random walk model for option pricing
\cite{CoxRossRubinstein1979}), insurance (Poisson random sums as
total claim size in dynamic collective risk models
\cite{Cramer1930}, binomial random sums as total claim size in
static portfolio risk models, geometric sums in the
Pollaczek--Khinchin--Beekman representation of the ruin probability
within the framework of the classical risk process
\cite{Kalashnikov1997}), reliability theory for modeling rare events
\cite{Kalashnikov1997}. It is now a tradition to admit that the
distributions of elementary jumps of these random walks may have
very heavy tails. The problems considered in the present paper
correspond to the situation where the tails may be as heavy as
possible for the normal approximation to be still adequate. A very
important (if not crucial) argument in favor of consideration of
approximations to the concentration functions and the corresponding
bounds for their accuracy when the variance exists is that in
financial mathematics the variance (or the square root of it) is
often used as a synonym of the {\it volatility}. At the same time,
from the Chebyshev inequality it is easy to obtain the inequality
$$
{\sf D}\xi\ge{\textstyle\frac14}\sup_{z\ge 0}z^2[1-Q_{\xi}(z)]
$$
which relates the concentration function $Q_{\xi}(z)$ of a random
variable $\xi$ with its variance ${\sf D}\xi$. Therefore, in
financial applications, the concentration function of (logarithmic)
increments of a stock price or some other financial index can be
used as a considerably more informative characteristic of the
volatility.

The paper is organized as follows. Section 1 contains main
definitions, preliminary information and auxiliary results. In
Section 2 we present some estimates for the deviation of the
concentration functions of sums of independent random variables with
finite variances from the folded normal distribution function
without any assumptions concerning the existence of the moments of
summands of higher orders. The obtained results are extended to
Poisson-binomial and binomial random sums in Section 3. The case of
the Poisson random sums is considered in Section 4. Under the same
assumptions, the bounds are obtained for the approximation of the
concentration functions of general mixed Poisson random sums by the
corresponding limit distributions in Section 5. As corollaries,
bounds are obtained for the accuracy of approximation of the
concentration functions of geometric, negative binomial and Sichel
random sums by the exponential, the folded variance gamma and the
folded Student distribution in Sections 6, 7 and 8, respectively.
Numerical estimates of all the constants involved are written out
explicitly.

\renewcommand{\thesection}{\arabic{section}}

\setcounter{section}{0}

\section{Preliminary information and auxiliary results}


\smallskip

{\sc Lemma 1.} {\it Let $\eta$ and $\xi$ be two random variables
such that
$$
\sup_x|{\sf P}(\eta<x)-{\sf P}(\xi<x)|\le\delta,\eqno(2)
$$
where $\delta>0$. Then
$$
\sup_{z\ge 0}|Q_{\xi}(z)-Q_{\eta}(z)|\le 4\delta.
$$
}

\smallskip

For the {\sc proof} see \cite{Galieva2}.

\smallskip

{\sc Remark 1.} If instead of (1) the ``concentration function''
$\widetilde Q_{\xi}(z)$ of a random variable $\xi$ is defined as
$$
\widetilde Q_{\xi}(z) = \sup_{x\in\r}{\sf P}(x \le \xi < x+z),\ \ \
\ \ z>0,\eqno(3)
$$
then the corresponding analog of Lemma 1 can be proved with the
twice less constant. Namely, if $\eta$ and $\xi$ are two random
variables such that (2) holds, then
$$
\sup_{z\ge 0}|\widetilde Q_{\xi}(z)-\widetilde Q_{\eta}(z)|\le
2\delta,
$$
see \cite{Galieva1}

\smallskip

Recall the definition of a unimodal distribution due to
A.~Ya.~Khinchin. A random variable $\xi$ is said to have the {\it
unimodal distribution}, if there exists a point $x_0$ such that the
distribution function $F_{\xi}(x)$ of the random variable $\xi$ is
convex for $x<x_0$ and the function $1-F_{\xi}(x)$ is convex for
$x>x_0$. Moreover, in this case the point $x_0$ is called the {\it
mode} of the random variable $\xi$. It is easy to see that any
unimodal distribution function is continuous everywhere, possibly,
except for its mode.

The following statement was given in the book \cite{Hengartner}
without proof.

\smallskip

{\sc Lemma 2.} {\it Let $\xi$ be a random variable with the
symmetric unimodal distribution. Then for $z>0$ we have}
$$
Q_{\xi}(z)={\sf P}\Big(|\xi|<\frac{z}{2}\Big).
$$

\smallskip

A rigorous {\sc proof} of this statement can be found in, say,
\cite{Galieva2}.

\smallskip

Let $X_1,X_2,\ldots $ be independent random variables with ${\sf E}
X_i=0$ and $0<{\sf E} X_i^2\equiv\sigma_i^2<\infty$, $i=1,2,\ldots $
For $n\in\N$ denote
$$
S_n=X_1+\ldots +X_n,\ \ \ \ B_n^2=\sigma_1^2+\ldots +\sigma_n^2.
$$
Let $\Phi(x)$ be the standard normal distribution function,
$$
\Phi(x)=\frac{1}{\sqrt{2\pi}}\int\limits_{-\infty}^{x}e^{-z^2/2}dz,\
\ \ \ x\in\R.
$$
The {\it folded normal distribution function} will be denoted as
$\Phi_0(x)$,
$$
\Phi_0(x)=\begin{cases}2\Phi(x)-1,&x\ge0,\vspace{0.5mm}\\
0,&x<0.\end{cases}
$$
It is easy to see that if $\zeta$ is a random variable with the
standard normal distribution function, then $\Phi_0(x)={\sf
P}(|\zeta|<x)$.

Denote
$$
\Delta_n=\sup_x|{\sf P}(S_n<xB_n)-\Phi(x)|.
$$
Let $\mathcal{G}$ be the class of real functions $g(x)$ of $x\in\R$
such that\vspace{-2mm}
\begin{itemize}
\item $g(x)$ is even;\vspace{-2mm}
\item $g(x)$ is nonnegative for all $x$ and $g(x)>0$ for $x>0$;\vspace{-2mm}
\item $g(x)$ and $x/g(x)$ do not decrease for $x>0$.\vspace{-1mm}
\end{itemize}

In 1963 M. Katz \cite{Katz1963} proved that, whatever
$g\in\mathcal{G}$ is, if the random variables $X_1,X_2,\ldots $ are
identically distributed and ${\sf E} X_1^2g(X_1)<\infty$, then there
exists a finite positive constant $C_1$ such that
$$
\Delta_n\le C_1\cdot\frac{{\sf E} X_1^2g(X_1)}{\sigma_1^2
g\big(\sigma_1\sqrt{n}\big)}.\eqno(4)
$$
In 1965 this result was generalized by V.\,V.\,Petrov
\cite{Petrov1965} to the case of non-identically distributed
summands (also see \cite{Petrov1972}): whatever $g\in\mathcal{G}$
is, if ${\sf E} X_i^2g(X_i)<\infty$, $i=1,\ldots ,n$, then there
exists a finite positive constant $C_2$ such that
$$
\Delta_n\le \frac{C_2}{B_n^2g(B_n)}\sum_{i=1}^n{\sf E}
X_i^2g(X_i).\eqno(5)
$$

Everywhere in what follows the symbol $\mathbb{I}(A)$ will denote
the indicator function of an event $A$. For
$\varepsilon\in(0,\infty)$ denote
$$
L_n(\varepsilon)=\frac{1}{B_n^2}\sum_{i=1}^n{\sf E}
X_i^2\mathbb{I}(|X_i|\ge \varepsilon B_n),\ \ \
M_n(\varepsilon)=\frac{1}{B_n^3}\sum_{i=1}^n{\sf E}
|X_i|^3\mathbb{I}(|X_i|< \varepsilon B_n).
$$
In 1966 L.\,V.~Osipov \cite{Osipov1966} proved that there exists a
finite positive absolute constant $C_3$ such that for any
$\varepsilon\in(0,\infty)$
$$
\Delta_n\le C_3\big[L_n(\varepsilon)+M_n(\varepsilon)\big]\eqno(6)
$$
(also see \cite{Petrov1987}, Chapt V, Sect.\,3, theorem 7). This
inequality is of special importance. Indeed, it is easy to see that
$$
M_n(\varepsilon)\le\frac{\varepsilon}{B_n^2}\sum_{i=1}^{n}{\sf E}
X_i^2\mathbb{I}(|X_i|< \varepsilon B_n)\le\varepsilon.
$$
Hence, from (6) it follows that for any $\varepsilon\in(0,\infty)$
$$
\Delta_n\le C_3\big(\varepsilon+L_n(\varepsilon)\big).\eqno(7)
$$
But, as is well known, the Lindeberg condition
$$
\lim_{n\to\infty}L_n(\varepsilon)=0\ \ \text{for any
}\varepsilon\in(0,\infty)
$$
is a {\it criterion} of convergence in the central limit theorem.
Therefore, in terminology proposed by
V.\,M.\,Zolotarev~\cite{Zolotarev1997}, bound (7) is {\it natural},
since it relates the convergence {\it criterion} with the
convergence {\it rate} and its heft-hand and right-hand sides
converge to zero or diverge simultaneously.

In 1968 inequality (6) in a somewhat more general form was re-proved
by W.\,Feller \cite{Feller1968}, who used the method of
characteristic functions to show that $C_3\le6$.

A special case of (6) is the inequality
$$
\Delta_n\le C_3'\big[L_n(1)+M_n(1)].\eqno(8)
$$
In the book \cite{Petrov1972} it was demonstrated that $C_3\le
2C_3'$.

For identically distributed summands inequality (8) takes the form
$$
\Delta_n\le \frac{C_4}{\sigma_1^2}{\sf
E}X_1^2\min\bigg\{1,\,\frac{|X_1|}{\sigma_1\sqrt{n}}\bigg\}.\eqno(9)
$$

In the papers \cite{Paditz1980, Paditz1984} L.~Paditz showed that
the constant $C_4$ can be bounded as $C_4<4.77$. In 1986 in the
paper \cite{Paditz1986} he noted that with the account of lemma 12.2
from \cite{Bhat1982}, using the technique developed in
\cite{Paditz1980, Paditz1984}, the upper bound for $C_4$ can be
lowered to $C_4<3.51$.

In 1984 A.\,Barbour and P.\,Hall \cite{BarbourHall1984} proved
inequality (8) by Stein's method and, citing Feller's result
mentioned above, stated that the method they used gave only the
bound $C_3'\le 18$ (although the paper itself contains only the
proof of the bound $C_3'\le22$). In 2001 L.\,Chen and K.\,Shao
published the paper \cite{ChenShao2001} containing no references to
Paditz' papers \cite{Paditz1980, Paditz1984, Paditz1986} in which
the proved inequality (8) by Stein's method with the absolute
constant $C_3'=4.1$.

In 2011 V.\,Yu.\,Korolev and S.\,V.\,Popov  \cite{KP2011_3} showed
that there exist universal constants $C_1$ and $C_2$ which do not
depend on a particular form of $g\in\mathcal{G}$, such that
inequalities (4), (5), (8) and (9) are valid with $C_1=C_4\le
3.0466$ and $C_2=C_3'\le 3.1905$. This result was later improved by
the same authors in the papers \cite{KorolevPopov2011,
KorolevPopovDAN}, where it was shown that $C_1=C_2=C_4=C_3'\le
2.011$.

Moreover, in the paper \cite{KorolevPopovDAN} lower bounds were
established for the universal constants $C_1$ and $C_2$. Namely, let
$g$ be an arbitrary function from the class $\mathcal{G}$. Denote by
$\mathcal{H}_g$ the set of all random variables $X$ satisfying the
condition ${\sf E}X^2g(X)<\infty$. Denote
$$
C^*=\sup_{g\in\mathcal{G}}\sup_{{X_i\in\mathcal{H}_g,}\atop{i=1,\ldots
,n}}\frac{\Delta_n B_n^2g(B_n)}{\sum_{i=1}^n{\sf E} X_i^2g(X_i)}.
$$
It is easily seen that $C^*$ is the least possible value of the
absolute constant $C_2$ that provides the validity of inequality (5)
for all functions $g\in\mathcal{G}$ at once. In the paper
\cite{KorolevPopovDAN} it was proved that
$$
C^*\ge\sup_{z\ge 0}\Big|\frac{1}{1+z^2}-\Phi(-z)\Big|=0.54093\ldots
$$

In the recent paper \cite{DorofeevaKorolev2016} the results
mentioned above were improved and extended. First, it was shown that
one can take $C_3=C_3'$. Second, the upper bounds of the absolute
constants mentioned above were sharpened and it was shown that
$C_3\le 1.8627$. Third, these results were extended to
Poisson-binomial, binomial and Poisson random sums. Under the same
conditions, bounds were obtained for the accuracy of the
approximation of the distributions of mixed Poisson random sums by
the corresponding limit law. In particular, the bounds were
constructed for the accuracy of approximation of the distributions
of geometric, negative binomial and Poisson-inverse gamma (Sichel)
random sums by the Laplace, variance gamma and Student
distributions, respectively. All absolute constants were written out
explicitly. The main result of the paper \cite{DorofeevaKorolev2016}
can be formulated as follows.

\smallskip

{\sc Lemma 3.} {\it For any $n\in\N$ there holds the inequality}
$$
\Delta_n\le 1.8627\big[L_n(1)+M_n(1)].
$$

\section{Bounds for the concentration functions of non-random sums of
independent random variables with finite variances}

The main result of this section is the following

\smallskip

{\sc Theorem 1.} {\it For any $n\in\mathbb{N}$ and any
$\varepsilon>0$ there holds the inequality}
$$
\sup_{z\ge 0}\Big|Q_{S_n}(z)-\Phi_0\Big(\frac{z}{2B_n}\Big)\Big|\le
7.4508\big[L_n(\varepsilon)+M_n(\varepsilon)].
$$

\smallskip

{\sc Proof.} The desired assertion follows from Lemma 1 with
$\eta=S_n$, ${\sf P}(\xi<x)=\Phi(x/B_n)$ and Lemmas 2 and 3.

\smallskip

{\sc Corollary 1.} {\it For any $n\in\mathbb{N}$ and any
$\varepsilon>0$ there holds the inequality}
$$
\sup_{z\ge 0}\Big|Q_{S_n}(z)-\Phi_0\Big(\frac{z}{2B_n}\Big)\Big|\le
7.4508\big[\varepsilon+L_n(\varepsilon)].
$$

\smallskip

Actually Corollary 1 declares that as soon as the Lindeberg
condition holds, that is, the central limit theorem holds, the
concentration function of the sum of independent random variables
can be approximated by the folded normal distribution function with
the argument appropriately linearly transformed.

By the same reasoning as that used to prove Theorem 1, in which the
role of Lemma 3 is played by inequality 2 with the constant
sharpened in \cite{DorofeevaKorolev2016}, a result similar to
Theorem 1 can be obtained in terms of the function
$g\in\mathcal{G}$.

\smallskip

{\sc Theorem 2.} {\it Whatever a function $g\in\mathcal{G}$ is such
that ${\sf E} X_i^2g(X_i)<\infty$, $i\ge1$, for any $n\in\mathbb{N}$
there holds the inequality}
$$
\sup_{z\ge 0}\Big|Q_{S_n}(z)-\Phi_0\Big(\frac{z}{2B_n}\Big)\Big|\le
\frac{7.4508}{B_n^2g(B_n)}\sum_{i=1}^n{\sf E} X_i^2g(X_i).
$$

\section{Bounds for the concentration functions of Poisson-binomial and binomial random sums}

From this point on let $X_1,X_2,\ldots$ be independent {\it
identically distributed} random variables with ${\sf E} X_i=0$ and
$0<{\sf E} X_i^2\equiv \sigma^2<\infty$. Let $p_j\in(0,1]$ be
arbitrary numbers, $j=1,2,\ldots$. For $n\in\N$ denote
$\theta_n=p_1+\ldots+p_n$, $\vp_n=(p_1,\ldots,p_n)$. The
distribution of the random variable
$$
N_{n,\vp_n}=\xi_1+\ldots+\xi_n,
$$
where $\xi_1,\ldots,\xi_n$ are independent random variables such
that
$$
\xi_j=\begin{cases}1 & \text{ with probability } p_j,\cr
                    0 & \text{ with probability } 1-p_j,
      \end{cases},\ \ \ j=1,\ldots,n,
$$
is usually called Poisson-binomial distribution with parameters
$n;\vp_n$. Assume that for each $n\in\N$ the random variables
$N_{n,\vp_n},X_1,X_2,\ldots$ are jointly independent. The main
objects considered in this section are {\it Poisson-binomial random
sums} of the form
$$
S_{N_{n,\vp_n}}=X_1+\ldots+X_{N_{n,\vp_n}}.
$$
As this is so, if $N_{n,\vp_n}=0$, then we assume
$S_{N_{n,\vp_n}}=0$.

For $j\in\mathbb{N}$ introduce the random variables $\wX_j$ by
setting
$$
\wX_j=\begin{cases}X_j & \text{ with probability } p_j,\cr 0 &
\text{ with probability } 1-p_j.\end{cases}
$$
If the common distribution function of the random variables $X_j$ is
denoted $F(x)$ and the distribution function with a single unit jump
at zero is denoted $E_0(x)$, then, as is easily seen,
$$
{\sf P}\big(\wX_j<x\big)=p_jF(x)+(1-p_j)E_0(x),\ \ \ \
x\in\mathbb{R},\ j\in\mathbb{N}.
$$
It is obvious that ${\sf E}\wX_j=0$,
$$
{\sf D}\wX_j={\sf E}\wX_j^2=p_j\sigma^2.\eqno(10)
$$
In what follows the symbol $\eqd$ will denote coincidence of
distributions.

\smallskip

{\sc Lemma 4}. {\it For any $n\in\N$ and $p_j\in(0,1]$}
$$
S_{N_{n,\vp_n}}\eqd \wX_1+\ldots+\wX_n, \eqno(11)
$$
{\it where the random variables on the right-hand side of $(11)$ are
independent.}

\smallskip

For the {\sc proof} see \cite{DorofeevaKorolev2016}.

\smallskip

With the account of (10) and (11) it is easy to notice that
$$
{\sf D}S_{N_{n,\vp_n}}=\theta_n\sigma^2.
$$

\smallskip

{\sc Theorem 3.} {\it For any $n\in\N$ and $p_j\in(0,1]$, $j\in\N$,}
$$
\sup_{z\ge
0}\Big|Q_{S_{N_{n,\vp_n}}}(z)-\Phi_0\Big(\frac{z}{2\sigma\sqrt{\theta_n}}\Big)\Big|
\le\frac{7.4508}{\sigma^2}\,{\sf
E}X_1^2\min\bigg\{1,\,\frac{|X_1|}{\sigma\sqrt{\theta_n}}\bigg\}.
$$

\smallskip

{\sc Proof.} In \cite{DorofeevaKorolev2016} it was proved that
$$
\Delta_{n,\vp_n}\equiv \sup_x\big|{\sf
P}\big(S_{N_{n,\vp_n}}<x\sigma\sqrt{\theta_n}\big)-\Phi(x)\big|\le\frac{1.8627}{\sigma^2}\,{\sf
E}X_1^2\min\bigg\{1,\,\frac{|X_1|}{\sigma\sqrt{\theta_n}}\bigg\}.
$$
So, the desired assertion follows from Lemmas 1, 2 and (12).

\smallskip

{\sc Theorem 4.} {\it Whatever a function $g\in\mathcal{G}$ is such
that ${\sf E} X_1^2g(X_1)<\infty$, there holds the inequality}
$$
\sup_{z\ge
0}\Big|Q_{S_{N_{n,\vp_n}}}(z)-\Phi_0\Big(\frac{z}{2\sigma\sqrt{\theta_n}}\Big)\Big|
\le \frac{7.4508{\sf
E}X_1^2g(X_1)}{\sigma^2g(\sigma\sqrt{\theta_n})}.
$$

\smallskip

{\sc Proof.} This assertion follows from Lemmas 1, 2 and the
estimate
$$
\Delta_{n,\vp_n}\le 1.8627\frac{{\sf
E}X_1^2g(X_1)}{\sigma^2g(\sigma\sqrt{\theta_n})}
$$
proved in \cite{DorofeevaKorolev2016}.

\smallskip

In particular, if $p_1=p_2=\ldots=p$, then the Poisson-binomial
distribution with parameters $n\in\N$ and $\vp_n$ becomes the
classical binomial distribution with parameters $n$ and $p$:
$$
N_{n,\vp_n}\eqd N_{n,p},\ \ \ {\sf
P}(N_{n,p}=k)=C_n^kp^k(1-p)^{n-k},\ \ \ k=0,\ldots,n.
$$
In this case $\theta_n=np$, so that ${\sf D}S_{N_{n,p}}=np\sigma^2$.
Note that, as it was proved in \cite{DorofeevaKorolev2016}, if the
summands of the sums $S_n$ have identical distribution, then
inequalities (4) and (9) hold with the absolute constants equal to
1.8546. So, in the same way as Theorems 3 and 4 were proved, from
Lemmas 1 and 2 with the account of (4) and (9) we obtain the
following statements.

\smallskip

{\sc Corollary 2.} {\it For any $n\in\N$ and $p_j\in(0,1]$,
$j\in\N$,}
$$
\sup_{z\ge
0}\Big|Q_{S_{N_{n,p}}}(z)-\Phi_0\Big(\frac{z}{2\sigma\sqrt{np}}\Big)\Big|
\le\frac{7.4184}{\sigma^2}\,{\sf
E}X_1^2\min\bigg\{1,\,\frac{|X_1|}{\sigma\sqrt{np}}\bigg\}.
$$

\smallskip

{\sc Corollary 3.} {\it Whatever a function $g\in\mathcal{G}$ is
such that ${\sf E} X_1^2g(X_1)<\infty$, there holds the inequality}
$$
\sup_{z\ge
0}\Big|Q_{S_{N_{n,p}}}(z)-\Phi_0\Big(\frac{z}{2\sigma\sqrt{np}}\Big)\Big|
\le \frac{7.4184{\sf E}X_1^2g(X_1)}{\sigma^2g(\sigma\sqrt{np})}.
$$

\section{Bounds for the concentration functions of Poisson random sums}

In addition to the notation introduced above, let $\lambda>0$ and
$N_{\lambda}$ be the random variable with the Poisson distribution
with parameter $\lambda$:
$$
{\sf P}(N_{\lambda}=k)=e^{-\lambda}\frac{\lambda^k}{k!},\ \ \
k\in\mathbb{N}\cup \{0\}.
$$
Assume that for each $\lambda>0$ the random variables
$N_{\lambda},X_1,X_2,\ldots$ are jointly independent. Consider the
{\it Poisson random sum}
$$
S_{N_{\lambda}}=X_1+\ldots+X_{N_{\lambda}}.
$$
If $N_{\lambda}=0$, then we set $S_{N_{\lambda}}=0$. It is easy to
see that ${\sf E}S_{\lambda}=0$ and ${\sf
D}S_{\lambda}=\lambda\sigma^2$. The accuracy of the normal
approximation to the distributions of Poisson random sum was
considered by many authors, see the historical surveys in
\cite{KorolevShevtsova, ShevtsovaPoisson}. However, the analogs of
the Katz--Osipov-type inequalities (4) and (9) under relaxed moment
conditions were obtained only recently in
\cite{DorofeevaKorolev2016}. Namely, the following statement was
proved there. Denote
$$
\Delta_{\lambda}\equiv \sup_x\big|{\sf
P}\big(S_{\lambda}<x\sigma\sqrt{\lambda}\big)-\Phi(x)\big|.
$$

\smallskip

{\sc Lemma 5}. {\it For any $\lambda>0$ and any function
$g\in\mathcal{G}$ such that ${\sf E} X_1^2g(X_1)<\infty$ we have}
$$
\Delta_{\lambda}\le\frac{1.8546}{\sigma^2}\,{\sf
E}X_1^2\min\bigg\{1,\,\frac{|X_1|}{\sigma\sqrt{\lambda}}\bigg\},\eqno(12)
$$
$$
\Delta_{\lambda}\le 1.8546\frac{{\sf
E}X_1^2g(X_1)}{\sigma^2g(\sigma\sqrt{\lambda})}.\eqno(13)
$$

\smallskip

Using (12) and Lemmas 1 and 2 we obtain a bound for the accuracy of
the approximation of $Q_{S_{\lambda}}(z)$ by the folded normal
distribution function.

\smallskip

{\sc Theorem 5.} {\it For any $\lambda>0$}
$$
\sup_{z\ge0}\Big|Q_{S_{\lambda}}(z)-\Phi_0\Big(\frac{z}{2\sigma\sqrt{\lambda}}\Big)\Big|\le
\frac{7.4184}{\sigma^2}\,{\sf
E}X_1^2\min\bigg\{1,\,\frac{|X_1|}{\sigma\sqrt{\lambda}}\bigg\}.
$$

\smallskip

Whereas Lemmas 1, 2 and inequality (13) yield

\smallskip

{\sc Theorem 6.} {\it Whatever a function $g\in\mathcal{G}$ is such
that ${\sf E} X_1^2g(X_1)<\infty$, there holds the inequality}
$$
\sup_{z\ge0}\Big|Q_{S_{\lambda}}(z)-\Phi_0\Big(\frac{z}{2\sigma\sqrt{\lambda}}\Big)\Big|\le
7.4184\frac{{\sf E}X_1^2g(X_1)}{\sigma^2g(\sigma\sqrt{\lambda})}.
$$

\smallskip

The upper bound of the absolute constant used in Lemma 5 is {\it
uniform} over the class $\mathcal{G}$. In specific cases this bound
can be considerably sharpened. For example, it is obvious that
$g(x)\equiv|x|\in\mathcal{G}$. For such a function $g$ inequality
(13) takes the form of the classical Berry--Esseen inequality for
Poisson random sums, the best current upper bound for the absolute
constant in which is given in \cite{Shevtsova2014}:
$$
\Delta_{\lambda}\le 0.3031\frac{{\sf
E}|X_1|^3}{\sigma^3\sqrt{\lambda}},
$$
so, if the third moment of the summands exist, then instead of
Theorems 5 and 6 we can obtain the bound
$$
\sup_{z\ge0}\Big|Q_{S_{\lambda}}(z)-\Phi_0\Big(\frac{z}{2\sigma\sqrt{\lambda}}\Big)\Big|\le
1.2124\frac{{\sf E}|X_1|^3}{\sigma^3\sqrt{\lambda}}.
$$

\section{Bounds for the concentration functions of general mixed Poisson random sums}

In this section we extend the results of the preceding section to
the case where the random number of summands has the mixed Poisson
distribution. For convenience, in this case we introduce an
``infinitely large'' parameter $n\in\mathbb{N}$ and consider random
variables $N_n^{\star}$ such that for each $n\in\mathbb{N}$
$$
{\sf
P}(N_n^{\star}=k)=\int\limits_{0}^{\infty}e^{-\lambda}\frac{\lambda^k}{k!}d{\sf
P}(\Lambda_n<\lambda),\ \ \ k\in\mathbb{N}\cup\{0\},\eqno(14)
$$
for some positive random variable $\Lambda_n$. For simplicity $n$
may be assumed to be the scale parameter of the distribution of
$\Lambda_n$ so that $\Lambda_n=n\Lambda$ where $\Lambda$ is some
positive ``standard'' random variable in the sense, say, that ${\sf
E}\Lambda=1$ (if the latter exists).

Assume that for each $n\in\mathbb{N}$ the random variable
$N_n^{\star}$ is independent of the sequence $X_1,X_2,\ldots$. As
above, let $S_{N_n^{\star}}=X_1+\ldots+X_{N_n^{\star}}$ and if
$N_n^{\star}=0$, then $S_{N_n^{\star}}=0$.

From (14) it is easily seen that, if ${\sf E}\Lambda_n<\infty$, then
${\sf E}N_n^{\star}={\sf E}\Lambda_n$ so that ${\sf
D}S_n=\sigma^2{\sf E}\Lambda_n$.

Denote
$$
\Delta^{\star}_n\equiv\sup_x\bigg|{\sf
P}\big(S_{N_n^{\star}}<x\sigma\sqrt{{\sf
E}\Lambda_n}\big)-\int\limits_{0}^{\infty}
\Phi\Big(\frac{x}{\sqrt{\lambda}}\Big)d{\sf
P}\big(\Lambda_n<\lambda{\sf E}\Lambda_n\big)\bigg|.
$$
For $x\in\mathbb{R}$ introduce the function
$$
G_n(x)={\sf
E}\min\Big\{1,\,\frac{|x|}{\sigma\sqrt{\Lambda_n}}\Big\}={\sf
P}\Big(\Lambda_n<\frac{x^2}{\sigma^2}\Big)+\frac{|x|}{\sigma}{\sf
E}\frac{1}{\sqrt{\Lambda_n}}\mathbb{I}\Big(\Lambda_n\ge\frac{x^2}{\sigma^2}\Big).\eqno(15)
$$
The expectation in (15) exists since the random variable under the
expectation sign is bounded by 1. Of course, the particular form of
$G_n(x)$ depends on the particular form of the distribution of
$\Lambda_n$. In \cite{DorofeevaKorolev2016} the following statement
was proved.

\smallskip

{\sc Lemma 6.} {\it If ${\sf E}\Lambda_n<\infty$, then}
$$
\Delta^{\star}_n\le\frac{1.8546}{\sigma^2}{\sf
E}X_1^2G_n(X_1)=\frac{1.8546}{\sigma^2}{\sf
E}X_1^2\min\bigg\{1,\,\frac{|X_1|}{\sigma\sqrt{\Lambda_n}}\bigg\}=
$$
$$
=\frac{1.8546}{\sigma^2}\bigg[{\sf
E}X_1^2\mathbb{I}\big(|X_1|\ge\sigma\sqrt{\Lambda_n}\big)+{\sf
E}\frac{|X_1|^3}{\sigma\sqrt{\Lambda_n}}\mathbb{I}\big(|X_1|<\sigma\sqrt{\Lambda_n}\big)\bigg],
$$
{\it where the random variables $X_1$ and $\Lambda_n$ are assumed
independent.}

\smallskip

Taking into account that all scale mixtures of zero-mean normals are
symmetric and unimodal, using Lemmas 1, 2 and 6 we obtain the
following bound for the accuracy of the approximation of
$Q_{S_{N_n^{\star}}}(z)$ by the scale mixture of the folded normal
distribution function.

\smallskip

{\sc Theorem 7.} {\it If ${\sf E}\Lambda_n<\infty$, then}
$$
\sup_{z\ge0}\bigg|Q_{S_{N_n^{\star}}}(z)-\int\limits_{0}^{\infty}
\Phi_0\Big(\frac{z}{2\sqrt{\lambda}}\Big)d{\sf
P}\big(\Lambda_n<\lambda{\sf E}\Lambda_n\big)\bigg|\le
$$
$$
\le\frac{7.4184}{\sigma^2}{\sf
E}X_1^2G_n(X_1)=\frac{7.4184}{\sigma^2}{\sf
E}X_1^2\min\bigg\{1,\,\frac{|X_1|}{\sigma\sqrt{\Lambda_n}}\bigg\}=
$$
$$
=\frac{7.4184}{\sigma^2}\bigg[{\sf
E}X_1^2\mathbb{I}\big(|X_1|\ge\sigma\sqrt{\Lambda_n}\big)+{\sf
E}\frac{|X_1|^3}{\sigma\sqrt{\Lambda_n}}\mathbb{I}\big(|X_1|<\sigma\sqrt{\Lambda_n}\big)\bigg],
$$
{\it where the random variables $X_1$ and $\Lambda_n$ are assumed
independent.}

\smallskip

In the subsequent sections we will consider special cases where
$\Lambda_n$ has the exponential, gamma and inverse gamma
distributions.

\section{Bounds for the accuracy of approximation of the concentration functions
of geometric random sums by the exponential law}

In this section we consider sums of a random number of independent
random variables in which the number of summands $N_n^{\star}$ has
the geometric distribution with parameter $p=\frac{1}{1+n}$,
$n\in\mathbb{N}$:
$$
{\sf P}(N_n^{\star}=k)=\frac{1}{n+1}\Big(\frac{n}{n+1}\Big)^k,\ \ \
k\in\mathbb{N}\cup\{0\}.\eqno(16)
$$
As usual, we assume that for each $n\in\mathbb{N}$ the random
variables $N_n^{\star},X_1,X_2,\ldots$ are independent. We again use
the notation $S_{N_n^{\star}}=X_1+\ldots+X_{N_n^{\star}}$. If
$N_n^{\star}=0$, then we set $S_{N_n^{\star}}=0$. It is easy to see
that ${\sf E}N_n^{\star}=n$, ${\sf D}S_{N_n^{\star}}=n\sigma^2$.
Note that for any $k\in\mathbb{N}\cup\{0\}$
$$
{\sf P}(N_n^{\star}=k)=\frac{1}{n}\int\limits_{0}^{\infty}{\sf
P}(N_{\lambda}=k)\exp\Big\{-\frac{\lambda}{n}\Big\}d\lambda,
$$
where $N_{\lambda}$ is the random variable with the Poisson
distribution with parameter $\lambda$. This means that for
$N_n^{\star}$ representation (14) holds with $\Lambda_n$ being an
exponentially distributed random variable with parameter $\frac1n$.

In what follows we will use traditional notation
$$
\Gamma(\alpha,z)\equiv\int\limits_{z}^{\infty}y^{\alpha-1}e^{-y}dy,\
\ \ \gamma(\alpha,z)\equiv\int\limits_{0}^{z}y^{\alpha-1}e^{-y}dy,\
\ \text{ and } \ \
\Gamma(\alpha)\equiv\Gamma(\alpha,0)=\gamma(\alpha,\infty)
$$
for the upper incomplete gamma-function, the lower incomplete
gamma-function and the gamma-function itself, respectively, where
$\alpha>0$, $z>0$.

In the case under consideration
$$
\frac{1}{n}\int\limits_{0}^{\infty}\Phi_0\bigg(x\sqrt{\frac{n}{\lambda}}\bigg)\exp\Big\{-\frac{\lambda
}{n}\Big\}d\lambda=\int\limits_{0}^{\infty}\Phi_0\Big(\frac{x}{\sqrt{y}}\Big)e^{-y}dy=1-e^{\sqrt{2}x},\
\ \ x\ge0
$$
(see, e. g., lemma 12.7.1 in \cite{KorolevBeningShorgin2011}), that
is, the approximate distribution is exponential with parameter
$\sqrt{2}$.

At the same time, the function $G_n(x)$ (see (15)) has the form
$$
G_n(x)= 1-\exp\Big\{-\frac{x^2}{n\sigma^2}\Big\}+
\frac{|x|}{n\sigma}\int\limits_{x^2/\sigma^2}^{\infty}\frac{e^{-\lambda/n}}{\sqrt{\lambda}}d\lambda=
\gamma\Big(1,\frac{x^2}{n\sigma^2}\Big)+\frac{|x|}{\sigma\sqrt{n}}\Gamma\Big(\frac12,\,\frac{x^2}{n\sigma^2}\Big).
$$
So, from theorem 7 we obtain the following result.

\smallskip

{\sc Corollary 4.} {\it Let $N_n^{\star}$ have the geometric
distribution $(16)$. Then}
$$
\sup_{z\ge0}\Big|Q_{S_{N_n^{\star}}}(z)-1+\exp\Big\{-\frac{z}{2\sigma\sqrt{n}}\Big\}\Big|\le
\frac{7.4184}{\sigma^2}\bigg\{{\sf E}\Big[X_1^2
\gamma\Big(1,\frac{X_1^2}{n\sigma^2}\Big)\Big]+\frac{1}{\sigma\sqrt{n}}{\sf
E}\Big[|X_1|^3\Gamma\Big(\frac12,\,\frac{X_1^2}{n\sigma^2}\Big)\Big]\bigg\}.
$$

\section{Bounds for the accuracy of approximation of the concentration functions
of negative binomial random sums by the folded variance-gamma
distribution}

The case more general than that considered in the preceding section
is the case of negative binomial random sums.

Let $r>0$ be an arbitrary number. Assume that representation (14)
holds with $\Lambda_n$ being a gamma-distributed random variable
with the density
$$
p(\lambda)=\frac{\lambda^{r-1}e^{-\lambda/n}}{n^r\Gamma(r)}\ \
\lambda>0.
$$
Then the random variable $N_n^{\star}$ has the negative binomial
distribution with parameters $r$ and $\frac{1}{n+1}$:
$$
{\sf
P}(N_n^{\star}=k)=\frac{1}{n^r\Gamma(r)}\int\limits_{0}^{\infty}e^{-\lambda}\frac{\lambda^k}{k!}\lambda^{r-1}e^{-\lambda/n}d\lambda=
\frac{\Gamma(r+k)}{\Gamma(r)\,k!}\Big(\frac{1}{1+n}\Big)^r\Big(\frac{n}{1+n}\Big)^k,\
\ \ \ \ k\in\mathbb{N}\cup\{0\}.\eqno(17)
$$
Let
$$
\mathcal{V}^+_r(x)\equiv\frac{1}{\Gamma(r)}\int\limits_{0}^{\infty}\!\Phi_0\Big(\frac{x}{\sqrt{\lambda}}\Big)\lambda^{r-1}e^{-\lambda}d\lambda,\
\ \ x\in\mathbb{R},
$$
be the folded symmetric variance-gamma distribution with shape
parameter $r$ (see, e. g., \cite{MadanSeneta1990}).

In the case under consideration ${\sf E}N_n^{\star}={\sf
E}\Lambda_n=nr$ so that ${\sf D}S_{N_n^{\star}}=nr\sigma^2$ and for
any $x\in\mathbb{R}$
$$
\int\limits_{0}^{\infty}\!\Phi_0\Big(x\sqrt{\frac{{\sf
E}\Lambda_n}{\lambda}}\Big)d{\sf
P}(\Lambda_n<\lambda)=\frac{1}{n^r\Gamma(r)}\int\limits_{0}^{\infty}\!\Phi_0\Big(x\sqrt{\frac{nr}{\lambda}}\Big)
\lambda^{r-1}e^{-\lambda/n}d\lambda=
$$
$$
=\frac{1}{\Gamma(r)}\int\limits_{0}^{\infty}\!\Phi_0\Big(\frac{x\sqrt{r}}{\sqrt{\lambda}}\Big)\lambda^{r-1}e^{-\lambda}d\lambda\equiv
\mathcal{V}^+_r(x\sqrt{r}).
$$

Here the function $G_n(x)$ (see (15)) has the form
$$
G_n(x)=
\frac{1}{n^r\Gamma(r)}\int\limits_{0}^{x^2/\sigma^2}\lambda^{r-1}e^{-\lambda/n}d\lambda+
\frac{|x|}{\sigma
n^r\Gamma(r)}\int\limits_{x^2/\sigma^2}^{\infty}\lambda^{r-3/2}e^{-\lambda/n}d\lambda=
$$
$$
=\frac{1}{\Gamma(r)}\Big[\gamma\Big(r,\frac{x^2}{n\sigma^2}\Big)+\frac{|x|}{\sigma\sqrt{n}}\Gamma\Big(r-\frac12,\frac{x^2}{n\sigma^2}\Big)\Big].
$$

So, from theorem 7 we obtain the following result.

\smallskip

{\sc Corollary 5.} {\it Let $N_n^{\star}$ have the negative binomial
distribution $(17)$. Then}
$$
\sup_{z\ge0}\Big|Q_{S_{N_n^{\star}}}(z)-\mathcal{V}^+_r\Big(\frac{z}{2\sigma\sqrt{n}}\Big)\Big|\le
\frac{7.4184}{\sigma^2\Gamma(r)}\bigg\{{\sf
E}\Big[X_1^2\gamma\Big(r,\frac{X_1^2}{n\sigma^2}\Big)\Big]+\frac{1}{\sigma\sqrt{n}}{\sf
E}\Big[|X_1|^3\Gamma\Big(r-\frac12,\frac{X_1^2}{n\sigma^2}\Big)\Big]\bigg\}.
$$

\section{Bounds for the accuracy of approximation of the concentration functions
of Poisson-inverse gamma random sums by the folded Student
distribution}

Let $r>1$ be an arbitrary number. Assume that representation (14)
holds with $\Lambda_n$ being an inverse-gamma-distributed random
variable with parameters $\frac{r}{2}$ and $\frac{n}{2}$ having the
density
$$
p(\lambda)=\frac{n^{r/2}\lambda^{-r/2-1}}{2^{r/2}\Gamma(\frac{r}{2})}\exp\Big\{-\frac{n}{2\lambda}\Big\},\
\  \ \lambda>0.
$$
Then the random variable $N_n^{\star}$ has the so-called
Poisson-inverse gamma distribution:
$$
{\sf
P}(N_n^{\star}=k)=\frac{n^{r/2}}{2^{r/2}\Gamma(\frac{r}{2})}\int\limits_{0}^{\infty}e^{-\lambda}\frac{\lambda^k}{k!}\lambda^{-r/2-1}
\exp\Big\{-\frac{n}{2\lambda}\Big\}d\lambda,\ \ \ \ \
k\in\mathbb{N}\cup\{0\},\eqno(35)
$$
which is a special case of the so-called Sichel distribution see, e.
g., \cite{Sichel1971, Willmot1993}. In this case
$$
{\sf E}\Lambda_n=\frac{n}{r-2}
$$
so that
$$
{\sf D}S_n^{\star}=\frac{n\sigma^2}{r-2}.
$$
Nevertheless, we will normalize random sums not by their mean square
deviations, but by slightly different and asymptotically equivalent
quantities $\sigma\sqrt{n/r}$.

As is known, if $\Lambda_n$ has the inverse gamma distribution with
parameters $\frac{r}{2}$ and $\frac{n}{2}$, then $\Lambda_n^{-1}$
has the gamma distribution with the same parameters. Therefore, we
have
$$
\frac{n^{r/2}}{\Gamma(\frac{r}{2})}\int\limits_{0}^{\infty}\Phi_0
\Big(x\sqrt{\frac{n}{r\lambda}}\Big)\lambda^{-r/2-1}\exp\Big\{-\frac{n}{2\lambda}\Big\}d\lambda=
\frac{n^{r/2}}{\Gamma(\frac{r}{2})}\int\limits_{0}^{\infty}\Phi_0
\Big(x\sqrt{\frac{n\lambda}{r}}\Big)\lambda^{r/2-1}\exp\Big\{-\frac{n\lambda}{2}\Big\}d\lambda=
$$
$$
=\frac{1}{2^{r/2}\Gamma(\frac{r}{2})}\int\limits_{0}^{\infty}\Phi\Big(x\sqrt{\frac{\lambda}{r}}\Big)\lambda^{r/2-1}e^{-\lambda/2}d\lambda=
\mathcal{T}^+_r(x),\ \ \ x\in\mathbb{R},
$$
where $\mathcal{T}^+_r(x)$ is the folded Student distribution
function with parameter $r$ ($r$ ``degrees of freedom'')
corresponding to the density
$$
t^+_r(x)=\frac{2\Gamma(\frac{r+1}{2})}{\sqrt{\pi
r}\Gamma(\frac{r}{2})}\Big(1+\frac{x^2}{r}\Big)^{-(r+1)/2},\ \ \
x\ge0,
$$
see, e. g., \cite{BeningKorolev2004}.

In this case the function $G_n(x)$ (see (15)) has the form
$$
G_n(x)={\sf
P}\Big(\Lambda_n^{-1}>\frac{\sigma^2}{x^2}\Big)+\frac{|x|}{\sigma}{\sf
E}\sqrt{\Lambda_n^{-1}}\mathbb{I}\mathbb{I}\Big(\Lambda_n^{-1}\le\frac{\sigma^2}{x^2}\Big)=
$$
$$
=\frac{n^{r/2}}{2^{r/2}\Gamma(\frac{r}{2})}\int\limits_{\sigma^2/x^2}^{\infty}\lambda^{r/2-1}e^{-n\lambda/2}d\lambda+
 \frac{|x|n^{r/2}}{2^{r/2}\sigma\Gamma(\frac{r}{2})}\int\limits_{0}^{\sigma^2/x^2}\lambda^{(r-1)/2}e^{-n\lambda/2}d\lambda=
$$
$$
=
\frac{1}{\Gamma(\frac{r}{2})}\Big[\Gamma\Big(\frac{r}{2},\frac{n\sigma^2}{2x^2}\Big)+
\frac{|x|}{\sigma}\sqrt{\frac{n}{2}}\,\gamma\Big(\frac{r+1}{2},\frac{n\sigma^2}{2x^2}\Big)\Big],
$$
where $\gamma(\,\cdot\,,\,\cdot\,)$ and
$\Gamma(\,\cdot\,,\,\cdot\,)$ are the lower and upper incomplete
gamma-functions, respectively. So, from theorem 7 we obtain the
following result.

\smallskip

{\sc Corollary 6.} {\it Let $N_n^{\star}$ have the Poisson-inverse
gamma distribution $(35)$. Then}
$$
\sup_{z\ge0}\Big|Q_{S_{N_n^{\star}}}(z)-\mathcal{T}^+_r\Big(\frac{z\sqrt{r}}{2\sigma\sqrt{n}}\Big)\Big|\le
\frac{7.4184}{\sigma^2\Gamma(\frac{r}{2})}\bigg\{{\sf
E}\Big[X_1^2\Gamma\Big(\frac{r}{2},\frac{n\sigma^2}{2X_1^2}\Big)\Big]+\frac{1}{\sigma}
\sqrt{\frac{n}{2}}\,{\sf
E}\Big[|X_1|^3\gamma\Big(\frac{r+1}{2},\frac{n\sigma^2}{2X_1^2}\Big)\Big]\bigg\}.
$$

\smallskip

{\sc Remark 2.} In accordance with Remark 1, all the theorems and
corollaries proved in this paper for the concentration functions $Q$
defined by relation (1), remain valid for the ``concentration
functions'' $\widetilde Q$ defined by (3). Moreover, the absolute
constants in the corresponding inequalities for $\widetilde Q$ are
twice less than those in the theorems proved for $Q$.

\small

\renewcommand{\refname}{Список литературы}

\end{document}